\documentclass[fleqn]{mat01}
\usepackage{times,mathtimy,amssymb}
\begin{document}

\setcounter{page}{159}
\firstpage{159}

\font\zz=msam10 at 10pt
\def\Box{\mbox{\zz{\char'244}}}

\newcommand{\f}{\frac}
\newcommand{\s}{\sigma}
\newcommand{\la}{\lambda}
\newcommand{\Si}{\Sigma}
\newcommand{\iy}{\infty}
\newcommand{\what}{\widehat}
\newcommand{\lgra}{\longrightarrow}

\newtheorem{theo}{Theorem}
\renewcommand\thetheo{\arabic{section}.\arabic{theo}}
\newtheorem{theor}[theo]{\bf Theorem}
\newtheorem{lem}[theo]{Lemma}
\newtheorem{pot}[theo]{Proof of Theorem}
\newtheorem{propo}{\rm PROPOSITION}
\newtheorem{rema}[theo]{Remark}
\newtheorem{defn}[theo]{\rm DEFINITION}
\newtheorem{exam}{Example}
\newtheorem{coro}[theo]{\rm COROLLARY}
\def\conjecture{\trivlist\item[\hskip\labelsep{\it Conjecture.}]}
\def\exammp{\trivlist\item[\hskip\labelsep{\it Example.}]}

\renewcommand{\theequation}{\thesection\arabic{equation}}

\newcommand{\R}{\mathbb R}%
\newcommand{\C}{\mathbb C}%
\newcommand{\D}{\mathbb D}%
\newcommand{\Z}{\mathbb Z}%
\newcommand{\Q}{\mbox{$\mathbb Q$}}%
\newcommand{\N}{\mathbb N}%

\newcommand{\cK }{\mbox{$ \mathcal K $}}%
\newcommand{\cL }{\mbox{$ \mathcal L $}}%
\newcommand{\cG }{\mbox{$ \mathcal G $}}%
\newcommand{\cN }{\mbox{$ \mathcal N $}}%
\newcommand{\cP }{\mbox{$ \mathcal P $}}%
\newcommand{\cR }{\mbox{$ \mathcal R $}}%
\newcommand{\sdp}{\mbox{\rule{. 1mm}{2mm}$\! \times$}}%

\newcommand{\Up }{\mbox{$ \upsilon $}}%
\newcommand{\aaa}{\mbox{$\{A_n\}^\infty_{n=1} \subseteq {\mathcal A}$}}%

\newcommand{\sss}{\mbox{$\displaystyle\sum^\infty_{n=1}$}}%
\newcommand{\ccc}{\mbox{$\displaystyle\coprod^\infty_{n=1}$}}%

\def\e{\mbox{\rm{e}}}
\def\d{\mbox{\rm{d}}}

\font\xxxxx=tir at 7.6pt
\def\ee{\mbox{\xxxxx{e}}}

\title{Cowling--Price theorem and characterization of heat kernel on
symmetric spaces}

\markboth{Swagato K Ray and Rudra P Sarkar}{Cowling--Price theorem}

\author{SWAGATO K RAY$^{1}$ and RUDRA P SARKAR$^{2}$}

\address{$^{1}$Department of Mathematics, Indian Institute
of Technology, Kanpur 208 016, India\\
\noindent $^{2}$Stat-Math Division, Indian Statistical Institute, 
203 B.T. Road, Kolkata 700 108, India\\
\noindent E-mail: skray@iitk.ac.in; rudra@isical.ac.in}

\volume{114}

\mon{May}

\parts{2}

\Date{MS received 13 February 2004; revised 12 April 2004}

\begin{abstract}
We extend the uncertainty principle, the Cowling--Price theorem, on
non-compact Riemannian symmetric spaces $X$. We establish a
characterization of the heat kernel of the Laplace--Beltrami operator
on $X$ from integral estimates of the Cowling--Price type.
\end{abstract}


\maketitle

\section{Introduction}

Our starting point in this paper is the classical Hardy's theorem
\cite{H}: if a measurable function $f$ on $\R$ satisfies $|f(x)| \leq
C\e^{-ax^2}$, $|\hat f(x)|\leq C\e^{-bx^2}$, $x\in\R$ for positive $a$
and $b$ with $ab>\f 14$, then $f=0$ almost everywhere. But if $a \cdot b
= \f14$ then $f$ is a constant multiple of the Gauss kernel
$\e^{-a|x|^2}$. The first assertion of the theorem, thanks to a number
of articles in the recent past (e.g. \cite{CSS,EEKK1,Se,SS}), may now be
viewed as instances of a fairly general phenomenon in harmonic analysis
of Lie groups known as Hardy's uncertainty principle. At the same time,
several variants of the above decay conditions have been employed in the
studies by many mathematicians, under which the first assertion has been
proved. Notable among them are the integral estimates on $f$ and
$\what{f}$ introduced by Cowling and Price \cite{CP}. Due to intrinsic
difficulties, however, research remains incomplete in most of these
cases as to the second assertion of the theorem, that is, which are the
functions that satisfy the sharpest possible decay conditions and it is
this aspect that we take up for study in this paper.

In this article we consider the problem on a Riemannian symmetric space
$X$ of non-compact type. We realize $X$ as $G/K$, where $G$ is a
connected non-compact semisimple Lie group with finite center and $K$ is
a maximal compact subgroup of $G$. Thus a function on $X$ is a right
$K$-invariant function on $G$ and $\what{f}$ is the operator Fourier
transform on the space of representations of the group $G$. On these
spaces Narayanan and Ray have discovered \cite{NR1} that the correct
Hardy-like estimates involve the heat kernel of the Laplace--Beltrami
operator on $X$. In another paper they have also proved \cite{NR2} the
first assertion, i.e., the uncertainty theorem under the integral
estimates of Cowling and Price (C--P). As we take up here the case of
sharpest decay in the sense of C--P, our results reinforce their
findings, functions satisfying sharp C--P estimates indeed involve the
heat kernel. There is however a further element of intricacy, as we
arrive at a hierarchy of sharp decay conditions, instead of just one
sharpest decay. These are obtained by looking at the case $a\cdot b =
\f14$ and tempering the exponentials with suitable polynomials in the
estimates. Our main result in this direction is Theorem~4.1.

We go on to notice that since the matrix coefficients of the
principal series representations are explicitly known in the case
of rank 1 symmetric spaces, we can exploit the relation between
the group-Fourier transform and Jacobi transform on $\R$. This
allows us to vary the polynomial in the C--P estimates which now
characterize not only the bi-invariant heat kernel but also some
allied class of left $K$-finite functions on $X$ arising again as
the solutions of the heat equation of the Laplace--Beltrami
operator.

As is often the case with analysis of the non-compact symmetric
spaces,  there is a strong analogy with the Euclidean space $\R^n$
looked upon as the homogeneous space $M(n)/O(n)$, where $M(n)$ is
the Euclidean motion group and $O(n)$ is the orthogonal group.
Here the Laplace operator $\Delta_n$ generates all
$M(n)$-invariant differential operators on $\R^n$ and the plane
waves $\e^{-i\lambda x\omega}$ can be thought of as the basic
eigenfunctions of $\Delta_n$. We consider the Fourier transform of
a function $f$ in polar coordinates: $\what{f}(\lambda, \omega)
= (2\pi)^{-n/2}\int_{\R^n} \e^{-i\lambda x \cdot \omega} f(x)\d x$ for
$\lambda\in \R^{+}, \omega\in S^{n-1}$ and further expand
$\what{f} (\lambda, \cdot)$\vspace{-.05pc} in terms of the eigenfunctions of the
Laplace operator of $S^{n-1}$ -- we call the coefficients of this
expansion as Fourier coefficients of $\what{f}(\lambda, \cdot)$.
Then we can replace the integral estimate on $\what{f}$ by similar
estimates on the Fourier coefficients of the function
$\what{f}(\lambda, \cdot)$ on $S^{n-1}$\vspace{-.05pc} and obtain a different version
of the C--P type result for $\R^n$ characterizing the heat kernel
of the Laplace operator $\Delta_n$\break (Theorem~3.3).

Back on the symmetric space $X$, we recall that the symmetric
spaces $X=G/K$ of rank 1 are also non-compact two-point
homogeneous spaces (see \cite{He1}) like $\R^n$. We can go over to
the analogue of the polar coordinates through the $KAK$
decomposition of $G$. It follows that the space $K/M$ is
$S^{m_\gamma + m_{2\gamma}}$ where $M$ is the centralizer of $A$ in
$K$ and  $m_\gamma, m_{2\gamma}$ are the multiplicities of the
two positive roots $\gamma$ and $2\gamma$ respectively. We impose
conditions on the Fourier coefficients of the Helgason Fourier
transform $\tilde{f}(\lambda, \cdot)$ treating them as functions on
$K/M$ and get a version of the C--P result which corresponds to
that of the Euclidean case  mentioned above. We prefer to put our
results for the rank 1 case in this direction instead of reworking
the results of the general case. Our result in this case can be
thought of as vindicating Helgason's programme initiated in his
Paley--Wiener theorem\break (see \cite{He2}).

Hardy's theorem on semisimple Lie groups was first taken up by
Sitaram and Sundari \cite{SS}. It inspired many articles in
recent times including ours (see \cite{FS} for a survey of  these
results). After this work was finished, we came to know about a
recent paper by Bonami {\it et~al} \cite{BDJ}. They have obtained a
Cowling--Price type result on $\R^n$, as a corollary of a more
general theorem of Beurling which they  prove for $\R^n$ in this
paper. Our result on $\R^n$ is only slightly different from what
they obtain in \cite{BDJ}. However note that  Beurling's theorem
has not yet been proved  for symmetric spaces. We may also mention
here that the last part of our work is influenced by a recent
paper by Thangavelu \cite{Th}.

Throughout this article for a $p\in[1,\infty]$, $p'$ denotes its
conjugate,\vspace{-.05pc} i.e. $\f 1p+\f 1{p'}=1$. For two functions $f_1$ and
$f_2$, $f_1(x)\asymp f_2(x)$ means there exists two positive
constants $C,C'$ such that $C f_2(x)\le f_1(x)\le C'f_2(x)$. We
follow the practice of using $C,C'$ etc. to denote  constants
(real or complex)\vspace{-.05pc} whose value may change from one line to the
next. We use subscripts and superscripts of $C$ when needed to
indicate their dependence on parameters. Required preliminaries
are given at the beginning of the sections.

This paper has some overlap with a paper of Andersen \cite{And}.
(In fact the author is kind enough to refer to the pre-print
version of this paper.)

\pagebreak

\section{A complex analytic result}

In this section we shall prove a result of complex analysis which
will be useful throughout this paper.

\begin{lem}\label{lemma-1} Let ${\mathcal Q} = \{\rho \e^{i\psi }|\rho >
0, \psi \in (0, \frac{\pi}{2})\}$. Suppose a function $g$ is analytic on
${\mathcal Q}$ and continuous on the closure $\overline{\mathcal Q}$ of
${\mathcal Q}$. Also suppose that for $p\in [1,\infty)$ and for
constants $C > 0, a > 0$ and $m\geq 0${\rm ,} {\rm (i)} $|g(x + iy)|\leq
C\e^{ax^2}(1 + |z|)^m$ for $z = x + iy \in \overline{\mathcal Q}$ and
{\rm (ii)} $\int_0^{\infty} |g(x)|^{p} \d x < \iy$. Then for $\psi\in [0,
\frac{\pi}{2}]$ and $\sigma\in \R^+, \int_{\sigma}^{\sigma + 1}|g(\rho
\e^{i\psi })|\d \rho \leq C \max \{\e^a, (\sigma + 1)^{1/p}\}(\sigma +
2)^{2m}$.
\end{lem}

\begin{proof}
For $m = 0$ this is proved in \cite{CP}. For $m > 0$, we define $h(z) =
\f{g(z)}{(i + z)^{2m}}$ for $\mbox{Im}\,z > 0$. Then the lemma follows by
applying the case $m = 0$ to the function $h$.\hfill $\Box$
\end{proof}

\begin{lem}\label{lemma-2} Let  $g$ be an entire function on $\C^d$ such
that{\rm ,}
\begin{enumerate}
\renewcommand{\labelenumi}{\rm (\roman{enumi})}
\leftskip .1pc
\item $|g(z)|\leq C\e^{a| {\rm Re}\, z|^2}(1 + | {\rm Im}\, z|)^m$ for some
$m > 0, a > 0${\rm ,}\vspace{.3pc}

\item $\int_{\R^d} \frac{|g(x)|^q}{(1 + |x|)^s}| Q(x)|\d x < \iy${\rm ,}
for some $q\ge 1, s > 1$ and a polynomial $Q$ of degree $M$ in $d$
variables.
\end{enumerate}

Then $g$ is a polynomial. Moreover $\deg g\le \min\{m, \f{s - M -
d}{q}\}$ and if $s \le q + M + d${\rm ,} then $g$ is a constant.
\end{lem}

\begin{proof}
Once we prove that $g$ is a polynomial, it would be clear from (i) and
(ii) respectively that, $\deg g\le m$ and $\deg g < \f{s-M-d}{q}$ and
hence if $s < q + M + d$, then $g$ is a constant. 

We will first assume that $d=1$. 

Since for a scalar $\alpha, g-\alpha$ also satisfies the above
conditions, we may and will assume that $g(0)=0$. Consider the function
$h$ given by $g(z) = zh(z)$. Thus $h$ is an entire function which
satisfies (i). Define
\begin{equation*}
H(z) = \frac{h(z)}{(1-iz)^{s/q}}\quad \mbox{ for } z\in {\mathcal Q},
\end{equation*}
where ${\mathcal Q} = \{\rho \e^{i\psi }|\rho > 0, \psi \in (0,
\frac{\pi}{2})\}$. Then $H$ is analytic on ${\mathcal Q}$ and continuous
on $\overline{\mathcal Q}$. Since $|1-iz|^{s/q}\ge 1$ for $z\in
\overline{\mathcal Q}$, we have
\begin{equation*}
|H(z)| \leq C\e^{a({\rm Re}\, z)^2}(1 + |\mbox{Im }\, z|)^m
\quad \mbox{ for } z\in \overline{\mathcal Q}.
\end{equation*}
Now
\begin{align*}
\int_0^{\infty} |H(x)|^q \d x &= \int_0^{\infty} \frac{|h(x)|^q}{|1 -
ix|^{s}}\d x\\[.2pc]
&= \int_0^{\infty} \frac{|h(x)|^q}{(1 + x^2)^{s/2}}\d x < \infty,
\end{align*}
by continuity of the integrand and (ii). Thus $H$ satisfies the
conditions of Lemma~\ref{lemma-1} and hence
\begin{equation*}
\int_{\sigma}^{\sigma + 1} |H(\rho \e^{i\psi})|\d \rho \leq A \max
\{\e^a, (\sigma + 1)^{1/q}\} (\sigma + 2)^{2m}.
\end{equation*}
Since $|1-i\rho \e^{i\psi }|\leq (1+\rho )$, it follows that
\begin{equation*}
\int_{\sigma}^{\sigma + 1} |h(\rho \e^{i\psi})|\d \rho \leq A \max
\{\e^a, (\sigma + 1)^{1/q}\} (\sigma + 2)^{2m} (\sigma + 2)^{s/q}.
\end{equation*}
Considering the functions $H_1(z) = \overline{h
(\overline{z})}/(1-iz)^{s/q},$ $H_2(z) = \overline{h (-
\overline{z})}/(1-iz)^{s/q},$ $H_3(z) = h(-z)/(1-iz)^{s/q}$ for $z\in
{\mathcal Q}$ we get that for large $\sigma$
\begin{equation}
\int_{\sigma}^{\sigma + 1} |h(\rho \e^{i\psi} )|\d \rho \leq A'(\sigma +
2)^{2m} (\sigma + 2)^{(s+1)/q}, \mbox{ for all } \psi.\label{*}
\end{equation}
By Cauchy's integral formula
\begin{equation*}
|h^r(0)| \leq r!(2\pi)^{-1} \int_0^{2\pi} |h(\rho \e^{i\psi})| \rho^{-r}
\d\psi.
\end{equation*}
For large $\sigma$ integrating both sides with respect to $\rho$ in
$[\sigma, \sigma + 1]$, we get 
\begin{align*}
|h^r(0)| &\leq r!(2\pi )^{-1}\sigma^{-r} \int_{\sigma}^ {\sigma
+1}\int_0^{2\pi}|h(\rho \e^{i\psi})|\d\psi\, \d\rho\\[.2pc]
&\leq C r!(2\pi)^{-1} \sigma^{-r}(\sigma + 2)^{(s+1)/q + 2m}
\end{align*}
by (\ref{*}). Let $\sigma\rightarrow\infty$, then $h^r(0) = 0$ for all
$r$ after some stage. Thus $h$ is a polynomial and hence so is $g$.

The proof of the lemma now proceeds by induction. Assuming the statement
of the lemma to be true for $d = 1, 2, \dots, n-1$, we prove it for $d =
n$. To this end we write $z = (z_1, z_2, \dots, z_n)\in \C^n$ as $z =
(\tilde{z}, z_n)$ where $\tilde{z} = (z_1, z_2, \dots, z_{n-1}) \in
\C^{n-1}$. We notice that for $x\in \R^n, (1 + |x|)\le (1 +
|\tilde{x}|)(1 + |x_n|)$. Thus the hypothesis of the lemma implies
$\int_{\R^{n}} \frac{|g (\tilde{x}, x_n)|^q}{(1 + |\tilde{x}|)^s (1 +
|x_n|)^s}|Q (\tilde{x}, x_n)|\d\tilde{x}\, \d x_n < \iy$.

By Fubini's theorem, we have a subset $B_1\subset\R^{n-1}$ of full
measure such that for $\tilde{x}\in B_1, \int_{\R} \frac{|g(\tilde{x},
x_n)|^q}{(1 + |x_n|)^s}|Q(\tilde{x}, x_n)|\d x_n < \iy$. Further, since
$Q$ is a non-zero polynomial we may choose $B_1$ so that for each
$\tilde{x}\in B_1$, $Q(\tilde{x}, \cdot)$ is a non-zero polynomial of
one variable. Thus if $\tilde{x}\in B_1$, the entire function
$g_{\tilde{x}}$ of one variable $g_{\tilde{x}} = g(\tilde{x}, z)$
satisfies the condition (i) and (ii) of the lemma and consequently, is a
polynomial of one variable. Writing $g(\tilde{z}, z_n) =
\sum_{m=0}^\infty a_m (\tilde{z}) z_n^m, \tilde{z}\in \C^{n-1}, z_n\in
\C$, where each $a_m$ is entire function of $\tilde{z}$, what we have
just proved means that for $\tilde{x}\in B_1, a_m(\tilde{x})$ is zero
except for finitely many values of $m$ (depending on $\tilde{x}$). But
each $a_m$ being an entire function, if $a_m$ is not identically zero,
$a_m(\tilde{x})=0$ only if $\tilde{x}\in N_m$ where $N_m\subset
\R^{n-1}$ has measure zero. If for infinitely many $m$ say for
$m\in\{m_1 < m_2 < \cdots < m_i < \cdots\}, a_m \not\equiv 0$, we would
have $B_1 \subset \cup_{i=1}^\infty N_{m_i}$, a contradiction. This
shows that $g(\tilde{z}, z_n) = \sum_{m=0}^M a_m (\tilde{z}) z_n^m$ for
some positive integer $M$. A similar argument with the role of
$\tilde{z}$ and $z_n$ reversed and the induction hypothesis for $d = n -
1$ will give an upper bound on the degree of the monomials in $\tilde{z}
= (z_1, z_2, \dots, z_{n-1})$ occurring in $g$. Together, it is proved
that $g(\tilde{z}, z_n)$ is a polynomial separately in $\tilde{z}$ and
$z_n$ and hence a polynomial in $z = (\tilde{z}, z_n)\in \C^n$ (see
\cite{P}). The degree of $g$ is estimated at the beginning of the proof.
Thus the proof is complete. \hfill $\Box$
\end{proof}

\section{Euclidean spaces}

For $f\in L^1 (\R^n)$ let
\setcounter{equation}{0}
\begin{equation}
\what{f} (y) = (2\pi)^{-n/2} \int_{\R^n} f(x) \e^{-ix \cdot y} \d x,
\label{fouriertransform}
\end{equation}
where by $x \cdot y$ we mean the inner product of $x$ and $y$ in $\R^n$.
Let $p_t(x) = (4\pi t)^{-n/2} \e^{-|x|^2/4t}$, $t > 0$ be the heat kernel
associated with the Laplacian on $\R^n$. Then $\what{p}_t(y) =
(2\pi)^{-n/2} \e^{-t|y|^2}$.

We can have the Cowling--Price theorem \cite{CP} extended as follows:

\setcounter{theo}{0}
\begin{theor}[\!]\label{Eu-Th1} Let $f$ be a measurable function on
$\R^n$ which satisfies{\rm ,} for $p, q\in [1,\infty)$ and $s, t > 0$

{\rm (i)}\ $\int_{\R^n} \f{|f(x)|^p \ee^{\f{p}{4s}|x|^2}}{(1 + |x|)^k}
\d x < \infty$ and {\rm (ii)}\ $\int_{\R^n}\f{|\what{f} (\lambda)|^q
\ee^{qt|\lambda|^2}}{(1 + |\lambda|)^l} \d\lambda < \infty$ for any
$k\in (n, n + p]$ and $l\in (n, n + q]$. Then

\begin{enumerate}
\renewcommand\labelenumi{\rm (\alph{enumi})}
\leftskip 1.7pc
\item if $s < t$ then $f\equiv 0${\rm ,}\vspace{.2pc}

\item if $s = t$ then $f$ is a constant multiple of $p_t${\rm ,}\vspace{.2pc}

\item if $s > t$ then there exists infinitely many linearly
independent functions satisfying {\rm (i)} and {\rm (ii)}.
\end{enumerate}
\end{theor}

We do not give our proof of the theorem as it runs along the same lines
as the proof, later in this paper, of the corresponding theorem on
Riemannian symmetric spaces (Theorem~\ref{main-result}). Instead we
proceed to a variant of the Cowling--Price theorem where the usual
Fourier transform is replaced by spherical harmonic coefficients of
$\what{f} (\lambda, \cdot)$, where $\what{f}(\lambda, \omega) =
(2\pi)^{-n/2} \int_{\R^n} \e^{-i\lambda x \cdot \omega} f(x) \d x,
\lambda \in \R^+, \omega\in S^{n-1}$.

For a non-negative integer $m$, let ${\mathcal H}_m$ denote the space of
spherical harmonics of degree $m$ on $S^{n-1}$. For a fixed $S_m\in
{\mathcal H}_m$ and $f\in L^1(\R^n), n\ge 2$, the Fourier coefficients
of $f$ in the angular variable are defined by
\begin{equation}
f_m (|x|) = f_m(r) = \int_{S^{n-1}} f(rx') S_m(x') \d x'\quad \mbox{ for
almost every } r > 0. \label{m-radial-comp}
\end{equation}
In the Fourier domain we define
\begin{equation}
F_m (\lambda) = \lambda^{-m} \int_{S^{n-1}} \what{f} (\lambda, \omega)
S_m(\omega)\d\omega,\quad \lambda > 0.\label{F_m-defn}
\end{equation}

We also note that if $F_m=0$ for all $S_m\in {\mathcal H}_m$ for all
$m\in \N$, then by the uniqueness of the Fourier transform $f_m\equiv 0$
for all $S_m\in \mathcal H_m$, $m\in \N$ and so $f\equiv 0$.

\begin{lem} Let $f\in L^1(\R^n)$ and for a fixed $S\in \mathcal H_m,
F_m$ and $f_m$ are as defined above. Then
\begin{equation}
F_m (\lambda) = C \int_{\R^{n+2m}} f_m(|x|)|x|^{-m} \e^{-i\lambda x \cdot
\omega} \d x\quad \mbox { for any } \omega \in S^{n + 2m - 1}, \lambda >
0.\label{F_m-FT-of-radial}
\end{equation}
\end{lem}

\begin{proof}
We recall two results. If $f$ is radial on $\R^n, f(x) = f_0(|x|)$, then
$\what{f}$ is also radial and we have for $y\in \R^n, |y| = r > 0$:
\begin{equation}
\what{f}(y) = \what{f}(|y|) = \what{f}(r) = r^{-((n-2)/2)} \int_0^\infty
f_0 (s) J_{\f{n-2}{2}} (rs)s^{n/2}\d s, \label{SteinTh3.3}
\end{equation}
where $J_k, k \ge \f 12$ is the Bessel function (see \cite{SW},
pp.~154--155). We also have (\cite{He1}, p.~25 Lemma~3.6),
\begin{equation}
\int_{S^{n-1}} \e^{i\lambda rx'\omega} S_m(\omega) \d\omega = C_{n,m}
\f{J_{\f{n}{2} + m - 1}(\lambda r)}{(\lambda r)^{\f{n}{2} - 1}}
S_m(x').\label{classical-to-bessel}
\end{equation}
From (\ref{F_m-defn}) and the definition of the Fourier transform in
polar coordinates we have,
\begin{equation*}
F_m (\lambda) = (2\pi)^{-n/2} \lambda^{-m} \int_{\R^n} \left(\int_{S^{n-1}}
\e^{i\lambda x \cdot \omega} S_m (\omega) \d\omega\right) f(x) \d x.
\end{equation*}
Let $x = rx'$ where $x'\in S^{n-1}$ and $|x| = r$. By
(\ref{classical-to-bessel}), (\ref{m-radial-comp}) and
(\ref{SteinTh3.3}) we get,
\begin{align}
&F_m (\lambda)\nonumber\\[.2pc]
&\quad\ = C_{n,m} (2\pi)^{-n/2} \lambda^{-m} \int_0^\infty
\int_{S^{n-1}} \f{J_{\f{n}{2} + m - 1}(\lambda r)} {(\lambda
r)^{\f{n}{2}-1}} r^{n-1} S_m(x') f(rx') \d x' \d r\nonumber\\[.2pc]
&\quad\ = C_{n,m} (2\pi)^{-n/2} \lambda^{-m} \int_0^\infty f_m (r)
\f{J_{\f{n}{2} + m - 1} (\lambda r)}{(\lambda r)^{\f{n}{2}-1}} r^{n-1}
\d r\nonumber\\[.2pc]
&\quad\ = C_{n,m} (2\pi)^{-n/2} \lambda^{-(n + 2m - 2)/2} \int_0^\infty f_m (r)
r^{-m} J_{\f{n}{2} + m - 1} (\lambda r) r^{\f{n + 2m}2} \d r\nonumber\\[.2pc]
&\quad\ = C_{n,m} (2\pi)^{-n/2} \int_{\R^{n + 2m}} f_m (|x|) |x|^{-m} \e^{-i
\lambda \omega \cdot x} \d x.\label{bessel-tr}
\end{align}
This establishes the lemma.\hfill $\Box$
\end{proof}

\begin{theor}[\!] Let $f{\rm :}\ \R^n \lgra \C$ be a measurable
function such that for $p, q \in [1,\infty), s, t > 0$ and for each
non-negative integer $m$ and $S_m \in {\mathcal H}_m${\rm ,}\vspace{.2pc}

{\rm (i)}\ $\int_{\R^n} \f{|f(x)|^p \ee^{\f{p}{4s} |x|^2}}{(1 + |x|)^k}
\d x < \infty$ and {\rm (ii)}\ $\int_{\R^+} \f{|F_m (\lambda)|^q \ee^{qt
\lambda^2}}{(1 + \lambda)^l}\d\lambda < \infty$ for any $k\in (n, n +
p]$ and $l \in (1, 1 + q]$. Then

\begin{enumerate}
\renewcommand\labelenumi{\rm (\alph{enumi})}
\leftskip 1.7pc
\item if $s < t$ then $f\equiv 0${\rm ,}\vspace{.2pc}

\item if $s = t$ then $f$ is a constant multiple of $p_t${\rm ,}\vspace{.2pc}

\item if $s > t$ then there exist infinitely many linearly independent
functions satisfying {\rm (i)} and {\rm (ii)}.
\end{enumerate}
\end{theor}

\begin{proof}
It follows from (i) that $f\in L^1(\R^n)$ and let $f_m$ be defined as
above with respect to $S\in \mathcal H_m$.

Let $I = \int_0^\infty |f_m (r)|^p \f{\ee^{\f{p}{4s} r^2}}{(1 + r)^k}
r^{n-1} \d r$. Then by Holder's inequality and (i) we get
\begin{equation}
I \le C_m \int_{\R^n} \f{|f(x)|^p \e^{\f{p}{4s} |x|^2}}{(1 + |x|)^k} \d
x < \infty.\label{**}
\end{equation}

Using polar coordinates in (\ref{F_m-FT-of-radial}) and taking $\lambda
= u + iv \in \C$ we get
\begin{align*}
&|F_m (u + iv)|\\[.2pc]
&\quad\ \le \int_0^\infty \int_{S^{n + 2m - 1}} |f_m (r)| r^{n + 2m - 1}
r^{-m} \e^{|v|r} \d x'\, \d r\ (\mbox{as } |x'| = |\omega| = 1)\\[.2pc]
&\quad\ = C \int_0^\infty \f{|f_m (r)| \e^{\f1{4s} r^2} r^{\f{n -
1}{p}}}{(1 + r)^{k/p}} (1 + r)^{k/p}\ r^{\f{1 - n}{p}} r^{(n + m -
1)} \e^{-\f1{4s} r^2 + |v|r} \d r\\[.2pc]
&\quad\ \le C \left[\int_0^\infty \e^{(- \f1{4s} r^2 + |v| r)
p'}(1 + r)^{kp'/p}\ r^{(\f{1-n}{p})p'} r^{(n + m - 1) p'}\d r
\right]^{1/p'}\\[.2pc]
&\qquad\ \ \mbox{(by Holder's inequality and (\ref{**})})\\[.2pc]
&\quad\ = C \left[ \int_0^\infty \e^{(-\f{p'}{4s} (r - 2s |v|))^2} (1 +
r)^{kp'/p}\ r^{n - 1 + mp'} \d r \right]^{1/p'}
\e^{sv^2}\\[.2pc]
&\quad\ \le C \e^{sv^2} (1 + s|v|)^R \mbox{ for some non-negative } R.
\end{align*}
If $s\le t$ then $|F_m (u + iv)| \le C \e^{t v^2} (1 + t|v|)^R$. From
the existence of the integral (\ref{F_m-FT-of-radial}) for all $u + iv
\in \C$ estimate one can easily show that $F_m$ is an entire function.
Let $G(z) = F_m(z) \e^{tz^2}$. Then, $|G(z)| \le c \e^{t({\rm Re}\,z)^2} (1 +
t| \hbox{Im}\, z|)^R$ and $\int_\R \f{|G(x)|^q}{(1 + |x|)^l} \d x <
\infty$\break by (ii).

Applying Lemma \ref{lemma-2} for $n = 1$ on $G(z)$ we have $F_m
(\lambda) = C_m \e^{-t \lambda^2} P_m (\lambda)$ for $\lambda \in \R$,
where $P$ is a polynomial whose degree depends on $R$ and $l$. Since $1
< l \le 1 + q$, we see from (ii) that the polynomial is constant and
hence $F_m (\lambda) = C_m \e^{-t \lambda^2}$. Therefore $f_m (|x|) =
C_m |x|^m p_t (x)$. But if $s < t$, it follows that the integral $I$
does not exist unless $C_m = 0$. Since the argument applies to all $S_m
\in \mathcal H_m, m = 0, 1, \dots, f\equiv 0$. 

If $s = t$, then the integral in (i) is again infinite unless $C_m = 0$
whenever $m > 0$. Thus $f(x) = f_0 (x) = C_0 p_t(x)$ by the uniqueness
of Fourier transform.

When $s > t$, let $s > t_0 > t$ for some $t_0$. Let $H(x)$ be a solid
harmonic of degree $k \ge 1$, i.e. $H(r, \omega) = r^kS (w), r\in \R^+,
\omega \in S^{n-1}$, for some $S \in {\mathcal H}_k$. Then it is easy to
see that $H \cdot p_{t_0}$ satisfies conditions (i) and (ii) of the
theorem. This completes the proof.\hfill $\Box$
\end{proof}

\section{Symmetric spaces}
\setcounter{equation}{0}

Let $G$ be a connected non-compact semisimple Lie group with finite
center and let $K$ be a fixed maximal compact subgroup of $G$. Our
set-up for the rest of the paper is the Riemannian symmetric space
$X=G/K$ equipped with the intrinsic metric $d$. Whenever
convenient, we will treat a function $f$ on $X$ also as a right
$K$-invariant function on  $G$.

Let $G = KAN$ (resp. $\mathfrak g = \mathfrak k + \mathfrak a +
\mathfrak n$) be an Iwasawa decomposition of $G$ (resp. $\mathfrak g$)
where $\mathfrak {g, k, a, n}$ are the Lie algebras of $G, K, A$ and $N$
respectively. Let $\mathfrak a^*$ and $\mathfrak a^*_\C$ respectively be
the real dual of $\mathfrak a$ and its complexification. Let $P = MAN$
be the minimal parabolic subgroup corresponding to this Iwasawa
decomposition where $M$ is the centralizer of $A$ in $K$. Let $M'$ be
the normalizer of $A$ in $K$. Then $W = M'/M$ is the (restricted) Weyl
group of $(G, A)$ which acts on $\mathfrak a_\C$ and on its dual
$\mathfrak a^*_\C$. Let $\Si(\mathfrak g, \mathfrak a)$ be the set of
restricted roots, $\Si^+\subset \Si(\mathfrak g, \mathfrak a)$ be the
set of positive restricted roots which is chosen once for all and
$\Sigma^+_0 \subset \Sigma^+$ be the set of indivisible positive roots.
Let us denote the underlying set of simple roots by $\Delta_0$ and the
corresponding positive Weyl chamber in $\mathfrak a$ by $\mathfrak a^+$.
Then $G$ has the Cartan decomposition $G = K \overline{A^+} K$, where
$A^+ = \exp \mathfrak a^+$. Let \hbox{$\langle\,,\,\rangle$} be the Killing
form. Suppose the real rank of $G$ is $n$, i.e. $\dim \mathfrak a=n$ and
$\{H_1, H_2, \dots , H_n\}$ an orthonormal basis (with respect to
\hbox{$\langle\,,\,\rangle$}) of $\mathfrak a$. For $\lambda\in \mathfrak
a^*$, let $H_\lambda\in \mathfrak a$ correspond to $\lambda$ via
\hbox{$\langle\,,\, \rangle$}, i.e. $\lambda(H)=\langle H_\lambda, H\rangle$
and let $\lambda_j = \langle H_\lambda, H_j\rangle$. Then
$\lambda = (\lambda_1, \dots, \lambda_n)$ identifies $\mathfrak a^*$ with
$\R^n$ and hence $\mathfrak a^*_\C$ with $\C^n$.

Among the representations of $G$, those relevant for analysis on $X$ are
the so-called (minimal) spherical principal series representations
$\pi_\lambda, \lambda \in \mathfrak a^*_\C = \C^n$. The representation
$\pi_\lambda$ is induced from the representation $\xi_0 \otimes \exp
(\lambda) \otimes 1$ of $P = MAN$ (where $\xi_0$ and $1$ are the trivial
representations of $M$ and $A$ respectively and $\lambda\in \C^n$ acts
as a character of the vector subgroup $A$) and is unitary only if
$\lambda \in \mathfrak a^* = \R^n$. For all $\lambda\in \C^n$, we can
realize $\pi_\lambda$ on the Hilbert space $L^2(K/M)$ (compact picture)
where for $x\in G$,
\begin{align*}
(\pi_{\lambda}(x)u)(k) &= \e^{-(i \lambda + \rho_0) (H (x^{-1} k))}
u (K (x^{-1} k))\\[.2pc]
&= \e^{(i \lambda + \rho_0) (A (xK, kM))} u(K (x^{-1} k)),\quad k \in K,
u \in L^2 (K/M),
\end{align*}
where $\rho_0 = \f{1}{2} \Sigma_{\gamma \in \Sigma^+} m_\gamma \gamma,
m_\gamma$ being the multiplicity of the root $\gamma$.

Here for $g \in G, K(g) \in K, H(g) \in \mathfrak a$ and $N(g) \in N$
are the {\em parts} of $g$ in the Iwasawa decomposition $G = KAN$, i.e.
$g = K(g) \exp H(g) N(g)$. The vector-valued inner product $A (\cdot ,
\cdot)$ is defined as $A(gK, kM) = - H(g^{-1} k)$ (see \cite{He2}). We
can choose an orthonormal basis of $L^2 (K/M)$ consisting of $K$-finite
vectors. Let $u$ be an arbitrary element of this basis and let $e_0$ be
the $K$-fixed vector in it. From the action of $\pi_{\lambda} (x)$
defined above, it is clear that~for all $x\in G, \langle e_0,
\pi_{i\rho_0} (x) e_0 \rangle =1$ and when $u \neq \e_0, \langle u,
\pi_{i\rho_0} (x) e_0 \rangle =0$. For a suitable function $f$ on $X$,
let $\what{f} (\lambda) = \int_{G/K} f(x) \pi_{\lambda} (x) \d x$ denote
its Fourier transform with respect to $\pi_{\lambda}$. The $(u, e_0)$th
matrix coefficient of the operator $\what{f} (\lambda)$, denoted by
$\what{f}_{u} (\lambda)$ is\break given by
\begin{equation*}
\what{f}_{u} (\lambda) = \int_{X} f(x)\langle u, \pi_{\lambda} (x) e_0
\rangle \d x.
\end{equation*}
Let $\d k$ and $\d a$ respectively be the Haar measures on $K$ and $A$
and $\int_K \d k =1$. Let $d$ be the distance on $X$ induced by the
Riemannian metric on it. We define $\sigma (x) = \d (xK, o)$, where $o =
e K, x \in X$. Then $\sigma (\exp H) = |H| = \langle H, H\rangle^{1/2}$
for all $H \in \mathfrak a$. For $\lambda = (\lambda_1, \lambda_2,
\dots, \lambda_n) \in \C^n$, by $|\lambda|$ we shall mean its usual norm
while $|\lambda|_\R$ will stand for $\lambda_1^2 + \lambda_2^2 + \dots +
\lambda_n^2$.

The following estimate of Harish-Chandra (\cite{H-C1}, \S9) will be
useful for us:
\begin{equation}
\e^{-\rho_0 (H)} \leq \Xi (\exp H) \le C \e^{-\rho_0 (H)} (1 +
|H|)^{|\Sigma^+_0|}\quad \mbox{ for all } H \in \overline{\mathfrak
a^+_0}, \label{apriori1}
\end{equation}
where $\Xi(x)$ is $\phi_0$, i.e. the elementary spherical function with
parameter $0$ and $|\Sigma^+_0|$ is the cardinality of $\Sigma^+_0$.

For $a\in A$, by $\log a$ we shall mean an element in $\mathfrak a$ such
that $\exp \log a = a$. The Haar measure $\d x$ on $G$ can be normalized
so that $\d x = J(a) \d k_1\, \d a\, \d k_2$, where $J(a) = \Pi_{\gamma
\in \Si^+} (\e^{\gamma(\log a)} - \e^{-\gamma (\log a)})^{m_\gamma}$ is
the Jacobian of the Cartan decomposition of $G$. Clearly,
\begin{equation}
|J(a)| \le C\e^{2\rho_0 (\log a)}.\label{jacobian}
\end{equation}
For notational convenience we will from now on use $\lambda_\R$ and
$\lambda_I$ respectively for the real and the imaginary parts of
$\lambda \in \C$ (instead of the usual $\hbox{Re}\, \lambda$ and
$\hbox{Im}\, \lambda$).

Let $\delta\in \what{K}$ and let $d_\delta$ be the degree of $\delta$.
Let us assume that $u\in L^2(K/M)$ transforms according to $\delta$.
Then, from the well-known estimate of the elementary spherical function
$\phi_\lambda$ (see e.g. \cite{GV}, Prop.~4.6.1) and using the arguments
of Mili\v ci\'c (\cite{M}, p.~83) (see also \cite{Se}, 4.2) we have
\begin{align}
|\langle u, \pi_{\lambda}(a) e_0 \rangle| \le (d_\delta)^{1/2}
\phi_{\lambda_I} (x) \le C_\delta \e^{\lambda^+_I (\log a)} \Xi (a)\ \ 
\mbox{ for } a \in A^+ \mbox{ and } \lambda \in \mathfrak a^*_\C,
\label{basicfinal2}\nonumber\\
\end{align}
where $\lambda^+_I$ is the Weyl translate of $\lambda_I$ which is
dominant, i.e. belongs to the positive Weyl chamber.

We need the following estimate for the Plancherel measure $\mu (\lambda)
= |c (\lambda)|^{-2}$ ($c (\lambda)$ being the Harish-Chandra's
$c$-function) (see \cite{A3}, p.~394):
\begin{equation}\label{estimate-plancherel-measure}
|c(\lambda)|^{-2} \asymp \Pi_{\gamma \in \Si_0^{+}} \langle \lambda,
\gamma \rangle^2 (1 + |\langle \lambda, \gamma \rangle|)^{m_\gamma +
m_{2 \gamma} - 2} \quad \mbox{ for } \lambda \in \mathfrak a^*.
\end{equation}
Using the identification of $\mathfrak a^*_\C$ and $\mathfrak a^*$ with
$\C^n$ and $\R^n$ respectively, we will write an element $\lambda$ in
$\C^n$ (resp. $\R^n$) as $(\lambda_1, \lambda_2, \dots, \lambda_n)$
where $\lambda_i \in \C$ (resp. $\in \R$) for $i = 1, \dots, n$ and $V$
for the degree of the polynomial
\begin{equation}\label{degV}
\Pi_{\gamma \in \Si_0^{+}} \langle \lambda, \gamma \rangle^2 (1 +
\langle \lambda, \gamma \rangle)^{m_\gamma + m_{2 \gamma} - 2}.
\end{equation}
With this preparation we come to the main theorem of this section.

\setcounter{theo}{0}
\begin{theor}[\!]\label{main-result}
Let $f\in L^1(X)\cap L^2(X)$ and let $p, q\in [1,\infty)$. Suppose for
some $k, l, a, b \in (0, \infty)${\rm ,}
\begin{equation}
\int_{X} \f{|f(x) \Xi(x)^{\f{2}{p} - 1} \e^{a \s(x)^2}|^p}{(1 +
\s(x))^k} \d x < \iy \label{onfunction}
\end{equation}
and
\begin{equation}
\int_{\R^n} \f{\| \what{f} (\lambda) \|_2^q \e^{qb (\lambda_1^2 +
\lambda_2^2 + \dots + \lambda_n^2)}}{(1 + |\lambda|)^l} \mu (\lambda) \d
\lambda < \iy, \label{onft}
\end{equation}
where $\what{f} (\lambda) = \int_X f(x) \pi_\lambda (x) \d x$ is the
operator values of Fourier transform of $f$ at the spherical principal
series $\pi_\lambda$ and $\|\what{f} (\lambda)\|_2$ is the
Hilbert--Schimdt norm of $\what{f} (\lambda)$.

\begin{enumerate}
\renewcommand\labelenumi{\rm (\roman{enumi})}
\leftskip .1pc
\item If $a \cdot b = \f{1}{4}${\rm ,} then $\what{f}_{u} (\lambda) =
P_{u,b} (\lambda) \e^{-b(\lambda_1^2 + \dots + \lambda_n^2)}, \lambda
\in \C^n${\rm ,} for some polynomial $P_{u,b}$ with $\deg P_{u,b} <
\min \Big\{ \f{2 |\Sigma^+_0|}{p'} + \f{k}{p} + 1, \f{l - V - n}{q}\Big\}${\rm ,}
where $V$ is given by {\em (\ref{degV})}.\vspace{.2pc}

If also $l\le q + V + n${\rm ,} then $P_{u, b}$ is a constant.\vspace{.2pc}

\item If $a \cdot b > \f 14${\rm ,} then $f\equiv 0$.
\end{enumerate}
\end{theor}

\begin{proof}
Let us recall that $\what{f}_{u} (\lambda) = \int_{G} f(x) \langle u,
\pi_{\lambda}(x) e_0 \rangle \d x$ is the $(u, 0)$th matrix coefficient
of $\what{f} (\lambda)$ for $\lambda \in \mathfrak a^*_{\C}$, if the
integral exists.

We shall show that $\what{f}_{u} (\lambda)$ exists and is an entire
function in $\lambda\in \C^n$ and
\begin{equation}
|\what{f}_{u} (\lambda)| \le C \e^{b |\lambda_I|^2} (1 +
|\lambda_I|)^{k'}\quad \mbox { for all } \lambda \in \C^n,\ \mbox{ for
some } k' > 0.\label{step1}
\end{equation}
We rewrite the condition (\ref{onfunction}) as
\begin{equation}
\int_{A^+} \f{|f(h) \Xi(h)^{\f{2}{p} - 1} \e^{a |\log h|^2}|^p}{(1 +
|\log h|)^k} J(h) \d h < \iy \quad \mbox { for all } h \in A^+.
\label{onfunction_kak}
\end{equation}

Then using (\ref{basicfinal2}), for all $\lambda\in \C^n$, we have
\begin{align*}
|\what{f}_{u} (\lambda)| &= \left|\int_{A^+} f(h) \langle u,
\pi_{\lambda}(h) e_0 \rangle J(h) \d h \right|\\[.2pc]
&\le C \int_{A^+} \left| \f{f(h) \Xi(h)^{\f{2}{p} - 1} \e^{a| \log
h|^2}}{(1 + |\log h|)^{k/p}} \right| \e^{-a |\log h|^2}\\[.2pc]
&\quad\ \times (1 + |\log h|)^{k/p} \e^{\lambda^+_{I} (\log h)} \Xi (h)^{2 (1 -
\f{1}{p})} J(h) \d h.
\end{align*}

From (\ref{onfunction_kak}) and applying (\ref{apriori1}) and 
(\ref{jacobian}), we get for $k' = \f{2 |\Sigma^+_0|}{p'} + \f{k}{p}$,
\begin{align*}
|\what{f}_{u} (\lambda)| &\le C \cdot \bigg(\int_{A^+} \e^{-p'a|\log h|^2}
\e^{p' \lambda_I^+ (\log h)} \e^{-2 \rho_0 (\log h)}\\[.2pc]
&\quad\ \times (1 + |\log h|)^{k'p'} \e^{2 \rho_0 (\log h)} \d h \bigg)^{1/p'}\\[.2pc] 
&= C \cdot \left(\int_{A^+} \e^{-p' a |\log h|^2} \e^{p' \lambda^+_{I} (\log
h)} (1 + |\log h|)^{k' p'} \d h \right)^{1/p'}\\[.2pc]
&\le C \cdot \left( \int_{\mathfrak a} \e^{-p' a |H|^2} \e^{p' \lambda^+_{I}
(H)} \cdot (1 + |H|)^{k'p'} \d H \right)^{1/p'},
\end{align*}
where $H = \log h$ and $\d H$ is the Lebesgue measure on $\mathfrak a$.
We recall that $H_{\lambda^+_I}$ corresponds to $\lambda_I^+$ as
mentioned above so that $|\lambda_I^+| = |H_{\lambda^+_I}|$. Then,
\begin{align*}
&|\what{f}_{u} (\lambda)|\\
&\quad\ \le C \cdot \e^{\f{1}{4a}| H_{\lambda^+_I}|^2} \left(
\int_{\mathfrak a} \e^{-p' a \langle H - \f{1}{2a} H_{\lambda^+_I}, H -
\f{1}{2a} H_{\lambda^+_I} \rangle} \cdot (1 + |H|)^{k'p'}\d H
\right)^{1/p'}.
\end{align*}

From this using translation invariance of Lebesgue measure, for $\lambda
\in \mathfrak a^*_{\C} = \C^n$,
\begin{align*}
|\what{f}_{u} (\lambda)| &\le C \cdot \e^{\f{1}{4a} |H_{\lambda^+_I}|^2}
(1 + |H_{\lambda^+_I}|)^{k'} \left(\int_{\mathfrak a} \e^{-p' a |H|^2}
(1 + |H|)^{k'p'} \d H \right)^{1/p'}\\[.2pc]
&\le C \cdot (1 + |\lambda_I|)^{k'} \e^{\f{1}{4a} |\lambda_I|^2}
\int_{\mathfrak a} \e^{-a |H|^2} (1 + |H|)^{k'}\d H\ \  \mbox{ as }
|\lambda_I^+| = |\lambda_I|\\[.2pc]
&= C' \cdot (1 + |\lambda_I|)^{k'} \e^{b|\lambda_I|^2}\quad \left(
\mbox{ as } b = \f{1}{4a} \right).
\end{align*}

The analyticity of $\what{f}_u (\lambda)$ is a result of an usual
argument from Cauchy's integral formula and Fubini's theorem.

Let $g(\lambda) = \e^{b (\lambda_1^2 + \lambda_2^2 + \dots +
\lambda_n^2)} \what{f}_{u} (\lambda)$. Then $|g (\lambda)| \le C'
\e^{b|\lambda_\R|^2} (1 + |\lambda_I|)^{k'}$. Also from condition
(\ref{onft}), $\int_{\R^n} \f{|g(\lambda)|^q}{(1 + |\lambda|)^s} \mu
(\lambda)\d \lambda < \iy$.

It follows from (\ref{estimate-plancherel-measure}) that
$|c(\lambda)|^{-2}$ can be replaced by a polynomial. Hence by 
Lemma~\ref{lemma-2}, $g$ is a polynomial say $P_{u,b} (\lambda)$ and
$\deg P_{u, b} \le k', \deg P_{u,b} < \f{l - V - n}{q}$. Thus,
$\what{f}_{u} (\lambda) = P_{u,b}(\lambda) \cdot \e^{-b (\lambda_1^2 +
\lambda_2^2 + \dots + \lambda_n^2)}$, for all $\lambda \in \C^{n}$.

If $l \le q + V + n$, then clearly $P_{u, b}$ is a constant. This proves
(i).

If $a \cdot b > \f 14$, then we can choose positive constants $a_1, b_1$
such that $a > a_1 = \f 1{4b_1} > \f 1{4b}$. Then $f$ and $\what{f}$
also satisfy (\ref{onfunction}) and (\ref{onft}) with $a$ and $b$
replaced by $a_1$ and $b_1$ respectively. Therefore $\what{f}_{u}
(\lambda) = P_{u, b_1} (\lambda) \e^{-b_1 (\lambda_1^2 + \dots +
\lambda_n^2)}$. But then $\what{f}_u$ cannot satisfy (\ref{onft}) for
any $u$ unless $P_{u, b_1} \equiv 0$, which implies $f\equiv 0$.\hfill $\Box$
\end{proof}

The following related results come as immediate consequences of 
Theorem~4.1.

\begin{coro}{\em \cite{NR2}}$\left.\right.$\vspace{.5pc}

\noindent Let $f$ be a measurable function on $X$ and let $p, q\in
[1,\infty)$. Suppose for $a, b\in \R^+$ with $a \cdot b \ge \f14${\rm ,}
\begin{equation*}
\hskip -4pc {\rm (i)}\hskip 3pc \int_{X} |f(x) \Xi(x)^{\f{2}{p} - 1}
\e^{a \s(x)^2}|^p \d x < \iy, x\in X
\end{equation*}
and
\begin{equation*}
\hskip -4pc {\rm (ii)}\hskip 2.85pc \int_{\R^n} \|\what{f} (\lambda) \|^q
\e^{q b (\lambda_1^2 + \lambda_2^2 + \dots + \lambda_n^2)} \mu (\lambda)
\d \lambda < \iy, \lambda = (\lambda_1, \lambda_2, \dots, \lambda_n) \in
\R^n
\end{equation*}
then $f\equiv 0$.
\end{coro}

\begin{proof} If $a \cdot b = \f 14$, then clearly $f$ and $\what{f}$
satisfy the conditions of Theorem~4.1 with $l < q + V + n$ and hence
$\what{f}_u (\lambda) = C \e^{-\f{1}{4a} (\lambda_1^2 + \dots +
\lambda_n^2)}$. From (ii) it follows that $C = 0$.

The proof for the case $a \cdot b > \f 14$, proceeds as in 
Theorem~4.1(ii).\hfill $\Box$
\end{proof}

\begin{coro}{\em \cite{NR1,Shi}}$\left.\right.$\vspace{.5pc}

\noindent Let $f$ be a measurable function on $X$. Suppose for positive
constants $a, b$ and $C$ and for $r\ge 0$ with $a \cdot b = \f14${\rm ,}
\begin{equation*}
\hskip -4pc {\rm (i)}\hskip 3.1pc |f(x)| \le C \Xi (x) \e^{-a \s(x)^2} (1
+ \sigma(x))^r, x\in X
\end{equation*}
and
\begin{equation*}
\hskip -4pc {\rm (ii)}\hskip 2.85pc \| \what{f}(\lambda) \| \le C'
\e^{-b (\lambda_1^2 + \lambda_2^2 + \dots + \lambda_n^2)}, \lambda\in
\R^n
\end{equation*}
then $\what{f}_{0} (\lambda) = C \e^{-\f{1}{4a} (\lambda_1^2 + \dots +
\lambda_n^2)}$ for $\lambda = (\lambda_1, \lambda_2, \dots, \lambda_n)
\in \C^n$ and $\what{f}_u \equiv 0$ for any non-trivial $u$.
\end{coro}

\begin{proof} From (i) and (ii) we observe that $f$ and $\what{f}$
satisfy the estimates of Theorem~4.1 with $l \le q + V + n$. So
$\what{f}_{u} (\lambda) = C_u \e^{-\f{1}{4a} (\lambda_1^2 + \dots +
\lambda_n^2)}$ for any $\lambda \in \C^n$. But since for any non-trivial
$u$, the matrix coefficient function $\Phi_{i\rho}^{u, 0} (x) = \langle
\pi_{i \rho_0} (x) u, e_0 \rangle$ of $\pi_{i\rho_0}$ is identically
zero, it follows that $\what{f}_{u} (\lambda) \equiv 0$. The only
non-zero part is $\what{f}_0 (\lambda) = C_0 \e^{-\f{1}{4a} (\lambda_1^2
+ \dots + \lambda_n^2)}$ for\break $\lambda\in \C^n$.\hfill $\Box$
\end{proof}

\begin{coro}{\em \cite{SS}}$\left.\right.$\vspace{.5pc}

\noindent Let $f$ be a measurable function on $X$. Suppose for  positive
constants $a, b$ and $C$ with $a \cdot b > \f14${\rm ,}
\begin{align*}
\hskip -4pc {\rm (i)}\hskip 3.1pc &|f(x)| \le C \e^{-a\s(x)^2}, x\in X\\[.2pc]
\hskip -4pc {\rm (ii)}\hskip 2.85pc &\|\what{f} (\lambda)\| \le C' \e^{-b (\lambda_1^2 +
\lambda_2^2 + \dots + \lambda_n^2)}, \lambda\in \R^n
\end{align*}
then $f\equiv 0$.
\end{coro}

\begin{proof}
As $a > \f 1{4b}$ we can choose positive constants $a_1$ and $b_1$ such
that $a > a_1 = \f 1{4b_1} > \f 1{4b}$. Thus $f$ and $\what{f}$ satisfy
(i) and (ii) of the previous corollary with $a$ and $b$ replaced by
$a_1$ and $b_1$ respectively. Therefore it follows that $\what{f}_{0}
(\lambda) = C_1 \e^{-b_1 (\lambda_1^2 + \dots + \lambda_n^2)}$ and
$\what{f}_u \equiv 0$ for any non-trivial $u$. But as $b > b_1,
\what{f}_{0}$ as above cannot satisfy (ii) unless $C_1 = 0$.\hfill $\Box$
\end{proof}

\noindent{\it Sharpness of the estimates.}\ \ To complete the picture we
should consider the case $a \cdot b < \f 14$ and also show the
optimality of the factor $\Xi^{\f2p-1}$ considered in Theorem~4.1. We
provide an example below to show that if we substitute $\Xi^{\f{2}{p} -
\ell}$ for $\Xi^{\f{2}{p} - 1}$, with $\ell \in [0, 1)$ in
(\ref{onfunction}), then there are infinitely many linearly independent
functions, satisfying this modified estimate while their Fourier
transforms still satisfy (4.7) in Theorem~4.1. Also we will see that the
same example will show that if $a \cdot b < \f 14$ in the hypothesis of
Theorem~4.1, then there are infinitely many linearly independent
functions satisfying (4.6) and (4.7).\vspace{.5pc}

\begin{exammp} Let $G = SL_2(\C)$ and $K$ be its maximal compact subgroup
$SU(2)$ and $X = SL_2(\C)/SU(2)$. Then,
\begin{equation*}
A = \left\{ a_t = \left. \begin{pmatrix} \e^t &0\\[.2pc]
0 &\e^{-t}
\end{pmatrix} \right|\ t\in \R\right\}.
\end{equation*}
Let $\alpha$ be the unique element in $\Sigma^+$ given by
$\alpha (\log a_t) = 2t$ which occurs with multiplicity 2. Then 
$\sigma (a_t) = 2|t|$. Every $\lambda\in \C$ can  be identified with
an element in $\mathfrak a^*_\C$ by $\lambda = \lambda \alpha$. In
this identification, $\rho = 2$,  the unitary spherical principal
series representations are given by elements in $\R$, the
Plancherel measure $|c(\lambda)|^{-2} = |\lambda|^2$ and the
elementary spherical function $\phi_\lambda (a_t) = \sin(2 \lambda
t)/\lambda \sinh (2t)$ (see \cite{He1}, p.~432).

For a suitable function $f$ on $\R$, let $\tilde{f}$ be its
Euclidean Fourier transform and $C_c^\infty (\R)_{\rm even}$ be
the set of even functions in $C_c^\infty (\R)$.

We define a bi-invariant function $g$ on $G$ by prescribing its
spherical Fourier transform $\what{g} (\lambda) = \tilde{\psi} (2
\lambda) \what{h} (\lambda) P(2 \lambda)$ for $\lambda \in \R$, where
$\psi \in C_c^\infty (\R)_{\rm even}$ with support $[-\zeta, \zeta]$ for
some $\zeta > 0$, $\what{h} (\lambda) = \e^{-\lambda^2/4}$ for
$\lambda\in \R$ and $P$ is an even polynomial on $\R$. It follows from
the characterization of the bi-invariant functions in the Schwartz space
$S(G)$ (see \cite{GV}) that $g\in S(G)$, $g$ is bi-invariant and hence
can be thought of as a function on $X$.

Therefore Fourier inversion gives us
\begin{align*}
g(a_t) &= C \cdot \int_\R \tilde{\psi} (2 \lambda) \e^{-\lambda^2/4} P(2
\lambda) \f{\sin (2 \lambda t)}{\lambda \sinh (2t)} \lambda^2
\d\lambda\\[.2pc]
&= \f{C}{\sinh (2t)} \int_\R \tilde{\psi} (\lambda) \cdot
\e^{-\lambda^2/16} \lambda P(\lambda) \sin (\lambda t) \d\lambda\\[.2pc]
&= \f{C}{\sinh (2t)} (\psi_1 *_E h)(t),
\end{align*} 
where $h(t) = \e^{-4t^2}, \psi_1 \in C^\infty_c (\R)$ is the odd
function supported on $[-\zeta, \zeta]$ such that $\tilde{\psi_1}
(\lambda) = \lambda P(\lambda) \tilde{\psi} (\lambda)$ (i.e. $\psi_1$ is
certain derivative of $\psi$) and $*_E$ the Euclidean convolution.
Therefore for large $t$ and hence, choosing $C$ sufficiently large, for
all $t\in \R$,
\begin{align*}
|g(a_t)| &\le \f{C}{|\sinh 2t|} \cdot \e^{-4t^2} \e^{8\zeta t}\\[.2pc]
&\le C \cdot \e^{-\sigma (a_t)^2} \e^{- (1 - 4 \zeta) 2t}\\[.2pc]
&= C \cdot \e^{-\sigma (a_t)^2} \e^{- (1 - 4 \zeta) 2t}\\[.2pc]
&\le C \cdot \e^{-\sigma (a_t)^2} \Xi (a_t)^{(1 - 4 \zeta)}.
\end{align*}
Now if we choose $\zeta$ so that $\ell = (1 - 4 \zeta) > 0$, then for
all $x \in X, g$ satisfies
\begin{align}
|g(x)| \leq C \e^{-\sigma(x)^2} \Xi (x)^\ell (1 + \s(x))^M\quad \mbox{
for } \ell \in (0, 1) \mbox{ and for some } M > 0.\label{sharp}\nonumber\\
\end{align}
Its Fourier transform is $\what{g} (\lambda) = \tilde{\psi} (\lambda)
\what{h} (\lambda) P(\lambda), \lambda \in \R$. As $\tilde{\psi}$ is
bounded on $\R$, 
\begin{equation}
|\what{g} (\lambda)| \le C' \e^{-|\lambda|^2/4} (1 + |\lambda|)^N\quad
\mbox{ on }\R\; \mbox{ for some }N > 0. \label{sharpft}
\end{equation}

By (\ref{sharp}), (\ref{sharpft}), (4.1) and (4.2) this $g$ clearly
satisfies for $p, q \in [1,\infty)$ and $\ell \in (0, 1)$,
\begin{equation}
\int_{X} \f{|g (x) \Xi(x)^{\f{2}{p} - \ell} \e^{a \s(x)^2} |^p}{(1 +
\s(x))^k}\d x < \iy \label{on-function}
\end{equation}
and
\begin{equation}
\int_{\R} \f{|\what{g} (\lambda)|^q \e^{qb \lambda^2}}{(1 +
|\lambda|)^l} \mu (\lambda) \d \lambda < \iy, \label{on-ft}
\end{equation}
where $a = 1, b = \f14$ and $k > 3 + Mp, l > 3 + Nq$.

Since we can choose any $\psi \in C^\infty_c (\R)_{\rm even}$ and any
even polynomial $P(\lambda)$ to construct such a function $g$, we have
infinitely many linearly independent functions which satisfy
(\ref{on-function}) and (\ref{on-ft}).
\end{exammp}\vspace{.5pc}

\noindent {\it Case $a \cdot b < \f 14$.}\ \ Notice that from
(\ref{sharp}) and (\ref{sharpft}) it follows that there exist constants
$C_1, C_2 > 0$ for which the function $g$ constructed above satisfies
the estimates
\begin{equation}
|g(x)| \leq C_1 \e^{-\f{1}{2} \s(x)^2} \Xi(x) \label{sharp-alt}
\end{equation}
and
\begin{equation}
|\what{g} (\lambda)| \le C_2 \e^{-\f{1}{5} |\lambda|^2}\quad \mbox{ on }
\R. \label{ft-alt}
\end{equation}
Therefore
\begin{equation}
\int_{X} \f{|g (x) \Xi (x)^{\f{2}{p} - 1} \e^{a' \s(x)^2} |^p}{(1 +
\s(x))^k} \d x < \iy \label{onfunction-alt}
\end{equation}
and
\begin{equation}
\int_{\R^n} \f{|\what{g} (\lambda) |^q \e^{qb' \lambda^2}}{(1 +
|\lambda|)^l} \mu (\lambda) \d \lambda < \iy, \label{onft-alt}
\end{equation}
where $a' = \f {1}{2}$ and $b' = \f15$ and hence $a' \cdot b' < \f 14$.
In \cite{Sa} we have used this example for $p = q = \infty$ case.\vspace{.5pc}

\noindent{\it Characterization of the heat kernel.}\ \ The heat kernel
on $X$ is an analogue of the Gauss kernel $p_t$ on $\R^n$ (see \S3). Let
$\Delta$ be the Laplace--Beltrami operator of $X$. Then (see \cite{St},
Chapter~v), $T_t = \e^{t \Delta}, t > 0$ defines a semigroup
(heat-diffusion semigroup) of operators such that for any $\phi\in
C^\iy_c(X)$, $T_t\phi$ is a solution of $\Delta u = \partial u/\partial
t$ and $T_t \phi \lgra \phi$ a.e. as $t \lgra 0$. For every $t > 0, T_t$
is an integral operator with kernel $h_t$, i.e. for any $\phi \in
C^\iy_c (X), T_t \phi = \phi * h_t$. Then $h_t, t > 0$ are bi-invariant
functions and $h$ as a function of the variables $t \in \R^+$ and $x \in
G/K$ is in $C^\iy (G \times \R^+)$ satisfying the properties:

\begin{enumerate}
\renewcommand\labelenumi{(\roman{enumi})}
\leftskip .3pc
\item $\{h_t\hbox{:}\ t > 0\}$ form a semigroup under convolution $*$.
That is, $h_t * h_s = h_{t + s}$ for $t, s > 0$.\vspace{.2pc}

\item $h_t(x)$ is a fundamental solution of $\Delta u = \partial
u/\partial t$.\vspace{.2pc}

\item $h_t \in L^1(G) \cap L^\iy (G)$ for every $t > 0$.\vspace{.2pc}

\item $\int_{X} h_t(x) \d x =1$ for every $t > 0$.
\end{enumerate}

Thus we see that the heat kernel $h_t$ on $X$ retains all the nice
properties of the classical heat kernel. It is well-known that $h_t$ is
given by (see \cite{A4}):
\begin{equation}\label{heatkernel}
h_t (x) = \frac{1}{|W|} \int_{\mathfrak a^*} \e^{-t (|\lambda|_\R^2 +
|\rho_0 |_\R^2)} \phi_\lambda (x) \mu (\lambda) \d \lambda,
\end{equation}
where for $z = (z_1, \dots, z_n) \in \C^n, |z|_\R = z_1^2 + \dots +
z_n^2$, as defined earlier. That is, the spherical Fourier transform of
$h_t, \what{h_t} (\lambda)_{0, 0} = \e^{-t (|\lambda|_\R^2 +
|\rho_0|_\R^2)}$. It has been proved in \cite{A4} (Theorem~3.1(i)) that
for any $t_0 > 0$, there exists $C > 0$ such that
\begin{equation}\label{ankerestimate}
h_t (\exp H) \leq Ct^{-n/2} \e^{-|\rho_0|_\R^2 t - \langle \rho_0,
H \rangle - \f{|H|^2}{4t}} (1 + |H|^2)^{\f{d_X - n}{2}}
\end{equation}
for $t_0 \ge t > 0$ and $H \in \overline{\mathfrak a^+}$, where $d_X =
\dim X$.

Now the following elegant characterization of the heat kernel follows
from Theorem~4.1 but we present it in the form of a theorem to stress
the point.

\begin{theor}[\!] Let $f$ be a measurable function on $X$ such that for
some $t > 0$ and for $p, q \in [1, \infty)${\rm ,}
\begin{equation}
\int_{X} \f{|f(x) \Xi(x)^{\f{2}{p} - 1} \e^{\f{1}{4t} \s(x)^2} |^p}{(1 +
\s(x))^k} \d x < \iy \label{onfunctioncor}
\end{equation}
and
\begin{equation}
\int_{\R^n} \f{|\what{f} (\lambda) \e^{t |\lambda|^2}|^q}{(1 +
|\lambda|)^l} \mu (\lambda) \d \lambda < \iy,\label{onftcor}
\end{equation}
where $k > (d_X - n) p + 2|\Sigma_0^+ | + n, V + n < l \le q + V + n$
{\em ($V$ is as in (\ref{degV}))} and $\int_{X} f(x) \d x = 1$. Then $f
= h_t$.
\end{theor}

\begin{proof}
Let $a = 1/4t$ and $b = t$. Then it follows from Theorem~4.1 that
$\what{f}_{u} (\lambda) = C_u \cdot \e^{-t (\lambda_1^2 + \lambda_2^2 +
\dots + \lambda_n^2)} = C_u \cdot \e^{-t |\lambda|^2}$, because the
condition $l \le q + V + n$ forces the polynomial $P_{u, b}$ of
Theorem~4.1 to be a constant. But when $u\neq e_0$, then as noted in the
beginning of this section where we have described the representations,
that $\langle u, \pi_{i \rho_0} (x) e_0 \rangle = 0$ for all $x$ and
hence $\what{f}_{u} \equiv 0$. This implies that for any non-trivial
$\delta \in \what{K}_0, f_{\delta} \equiv 0$, where $f_{\delta} =
d_\delta \chi_\delta *_K\!\!f$ is the $\delta$th projection of $f$.
Thus $f$ is bi-invariant and its spherical Fourier transform $\what{f}
(\lambda)_{0, 0} = C_0 \cdot \e^{-t |\lambda|_\R^2}$. Again since,
$\phi_{\pm i \rho_0} (\cdot) \equiv 1, \int_{X} f(x) \d x = \int_{X}
f(x) \phi_{\pm i\rho_0} (x) \d x = \what{f}(\pm i \rho_0)_{0, 0}$.
Therefore the given initial condition reduces to $\what{f} (\pm i
\rho_0)_{0, 0} = 1$. From this, we have $\what{f} (\lambda)_{0, 0} =
\e^{-t (|\lambda|_\R^2 + |\rho_0|_\R^2)}$. This completes the proof in
view of (\ref{heatkernel}) above. It is clear from the estimates
(\ref{apriori1}) and (\ref{ankerestimate}) that $k$ and $l$ taken above
are good enough to\break accommodate $h_t$.\hfill $\Box$
\end{proof}

\section{Symmetric spaces of rank 1}

We shall revisit the Cowling--Price theorem from the point of view of
\S2 and relate it with the result obtained in the previous section. We
shall play around with the polynomials in the denominator of the
integrand of the C--P estimates. We will see that a larger class of
solutions of the heat equation can in fact be characterized using these
estimates, the usual heat kernel being one of them. The main technical
tool in this section will be the Jacobi functions. We shall crucially
use their relations with the Eisenstein integrals. Here the Jacobi
functions will take the role played by the Bessel functions in \S2.

Throughout this section the symmetric space $X$ is of rank 1. We will
continue to use the set-up and notation of the previous section, adapting
to this particular case. Here $\rho_0 = \f12 (m_{\gamma} + 2m_{2
\gamma})$. For $\lambda\in \C$, the function $x \lgra \e^{(i \lambda +
\rho_0) A(x, b)}$ is a common eigenfunction of invariant differential
operators on $X$. This motivates one to define Helgason--Fourier
transform of a function as a generalization of the Fourier transform in
polar\break coordinates:
\setcounter{equation}{0}
\begin{equation}
\tilde{f} (\lambda, b) = \int_X f(x) \e^{(-i \lambda + \rho_0) A(x, b)}
\d x, \label{Helgasontransform}
\end{equation}
where $\d x$ is the $G$-invariant measure on $X$. Let $\what{K}_0$ be
the set of equivalence class of irreducible unitary representations of
$K$ which are class 1 with respect to $M$, i.e. contains an $M$-fixed
vector -- it is also known that an $M$-fixed vector is unique up to a
multiple (see \cite{Kos}). Let $(\delta, V_\delta) \in \what{K}_0,
\delta$ different from the identity representation. Suppose $\{v_i | i =
1, \dots, d_\delta\}$ is an orthonormal basis of $V_\delta$ of which
$v_1$ is the $M$-fixed vector. Let $Y_{\delta, j} (kM) = \langle v_j,
\delta (k) v_1 \rangle, 1 \le j \le d_\delta$ and let $Y_0$ be the
$K$-fixed vector, which we have denoted by $e_0$ in the previous
section. Recall that $L^2 (K/M)$ is the carrier space of the spherical
principal series representations $\pi_\lambda$ in the compact picture
and \hbox{$\{Y_{\delta, j}\hbox{:}\ 1 \le j \le d_\delta, \delta \in
\what{K}_0\}$} is an ortho- normal basis for $L^2 (K/M)$ adapted to the
decomposition $L^2 (K/M) = \Si_{\delta \in \what{K}_0} V_\delta$ (see
\cite{He2}). As the space $K/M$ is $S^{m_\gamma + m_2 \gamma} =
S^{2\alpha + 1}$, this decomposition can be viewed as the spherical
harmonic decomposition and therefore $Y_{\delta, j}$'s can be considered
as the {\em spherical harmonics}.

For $\delta\in \what{K}_0, 1\le j\le d_\delta, \lambda\in \mathfrak
a^*_\C$ and $x\in X$, define,
\begin{equation}
\Phi_{\lambda, \delta}^j (x) = \int_K \e^{(i \lambda + \rho_0) A(x, kM)}
Y_{\delta, j} (kM) \d k.\label{Phi^j}
\end{equation}
Then $\Phi_{\lambda, \delta}^j(x)$ is a matrix coefficient of the
generalized spherical function (Eisenstein integral) $\int_K \e^{(i
\lambda + \rho_0) A(x, kM)} \delta(k) \d k$ (see \cite{He2}). Again,
$\Phi_{\lambda, \delta}^j (x) = \langle Y_{\delta, j}, \pi_{-\bar
\lambda} (x) Y_0 \rangle$, i.e. $\Phi_{\lambda, \delta}^j$ is a matrix
coefficient of the spherical principal series in the compact picture
(see \S4). It is well-known (see \cite{He2}) that they are
eigenfunctions of the Laplace--Beltrami operator $\Delta$ with
eigenvalues $-(\lambda^2 + \rho_0^2)$. When $\delta$ is trivial then
$\Phi_{\lambda, \delta}^1$ is obviously the elementary spherical
function $\phi_\lambda$.

The following result can be viewed as an analogue of
(\ref{classical-to-bessel}). For $\lambda\in \mathfrak a^*_\C$, $x = k
a_r K \in X$ and $1 \le j\le d_\delta$ (see \cite{He2}, p. 344)
\begin{equation}
\Phi^j_{\lambda, \delta} (x) = Y_{\delta, j} (kM) \Phi^1_{\lambda,
\delta} (a_r). \label{helgason-to-jacobi}
\end{equation}

Elements $\delta\in \what{K}_0$ can be parametrized by a pair of
integers $(p_\delta, q_\delta)$ so that $p_\delta \ge 0$ and $p_\delta
\pm q_\delta \in 2 \Z^+$ (see \cite{JW,Kos}). The trivial
representation in $\what{K}_0$ is parametrized by $(0, 0)$ in this 
set-up. Every $m > 0$ determines a subset $\what{K}_0 (m)$ of $\what{K}_0$
by $\what{K}_0 (m) = \{\delta \in \what{K}_0\hbox{:}\ p_\delta < m\}$.
This set is finite because, $p_\delta \ge |q_\delta|$. This
parametrization of $\what{K}_0$ will make a crucial appear- ance in our
results. We shall come back to that, after a digression.\vspace{.5pc}

\noindent {\it Jacobi functions.}\ \ At this point we give a quick
review of some preliminaries on the Jacobi functions. For a detailed
exposition the reader is referred to \cite{Koo}.

For $\alpha, \beta, \lambda \in \C, \alpha$ not a negative integer and
$r\in \R$, let $\phi_\lambda^{(\alpha, \beta)}(r)$ be the Jacobi
function of type $(\alpha, \beta)$ which is given in terms of the
hypergeometric function $_2\!F_1$ as
\begin{equation*}
\phi_\lambda^{(\alpha, \beta)}(r) = \,_2\!F_1 \left(\f{\alpha + \beta +
1 + i \lambda}{2}, \f{\alpha + \beta + 1 - i \lambda}{2}; \alpha + 1; -
(\sinh r)^2 \right).
\end{equation*}
Let
\begin{equation*}
\mathcal{L}_{\alpha, \beta} = \f{\d^2}{\d r^2} + ((2 \alpha + 1) \coth r
+ (2 \beta + 1) \tanh r) \f{\d}{\d r}
\end{equation*}
be the Jacobi Laplacian and let $\rho = \alpha + \beta + 1$. Then
$\phi_\lambda^{(\alpha, \beta)}(r)$ is the unique analytic solution of
the equation $\mathcal{L}_{\alpha, \beta} \phi = - (\lambda^2 + \rho^2)
\phi$, which is even and $\phi_\lambda^{(\alpha, \beta)}(0) = 1$. They
also satisfy respectively the following relation and the estimate: For
$r \in \R^+$ and $\lambda \in \C$,
\begin{equation}
\phi_\lambda^{\alpha, \beta} (r) = (\cosh r)^{-2 \beta}
\phi_\lambda^{\alpha, -\beta} (r) \label{jacobi-beta-to-minus-beta}
\end{equation}
and
\begin{equation}
|\phi_\lambda^{(\alpha, \beta)} (r)| \le C(1 + r) \e^{r (|\lambda_I| -
\rho)}. \label{estimate-jacobi-function}
\end{equation}

The associated Jacobi function $\phi_{\lambda, p, q}^{(\alpha, \beta)}$
for two extra parameters $p, q \in \Z$ is defined as
\begin{equation*}
\phi_{\lambda, p, q}^{(\alpha, \beta)} (r) = (\sinh r)^p (\cosh r)^q
\phi_{\lambda}^{(\alpha + p, \beta + q)} (r).
\end{equation*}
We also have the associated Jacobi operator,
\begin{align*}
\mathcal{L}_{\alpha, \beta, p, q} &= \f{\d^2}{\d r^2} + ((2 \alpha + 1)
\coth r + (2 \beta + 1) \tanh r) \f{\d}{\d r}\\[.2pc]
&\quad\ + \{-(2 \alpha + p) p (\sinh r)^{-2} + (2 \beta + q) q (\cosh
r)^{-2}\}.
\end{align*}
One can easily verify that $\phi_{\lambda, p, q}^{(\alpha, \beta)}$ is
again the unique solution of $\mathcal{L}_{\alpha, \beta, p, q} \phi = -
(\lambda^2 + \rho^2) \phi$, which is even and satisfies $\phi_{\lambda,
p, q}^{(\alpha, \beta)} (0) = 1$, in the case $p=0$. The proof of our
next theorem will involve finding a relation between the heat kernels of
the operators $\mathcal{L}_{\alpha + p,\, \beta + q}$ and
$\mathcal{L}_{\alpha, \beta, p, q}$.

Let
\begin{equation}\label{jacobi-delta}
\Delta_{\alpha, \beta} (r) = (2 \sinh r)^{2 \alpha + 1} (2 \cosh r)^{2
\beta + 1} = 4^\rho (\sinh r)^{2 \alpha + 1}(\cosh r)^{2 \beta + 1}.
\end{equation}
Then the (Fourier-) Jacobi transform of a suitable function $f$ on
$\R^+$ is defined by
\begin{equation}
J_{\alpha, \beta} (f) (\lambda) = \int_0^\infty f(r)
\phi_\lambda^{(\alpha, \beta)} (r) \Delta_{\alpha, \beta} (r) \d r,\,
\lambda\in \C.
\end{equation}

The {\em heat kernel} $h_t^{(\alpha, \beta)}$, associated to
$\mathcal{L}_{\alpha, \beta}$ is defined by
\begin{equation}\label{heat-kernel-jacobi-1}
\int_0^\infty h_t^{(\alpha, \beta)} (r) \phi_\lambda^{(\alpha, \beta)}
(r) \Delta_{\alpha, \beta} (r) \d r = \e^{-(\lambda^2 + \rho^2)t}.
\end{equation}
If $\alpha > - 1$ and $\alpha \pm \beta \geq - 1$, then we have {\em
inversion formula} for the Jacobi transform \cite{Koo} from which we
get
\begin{equation}
h^{(\alpha, \beta)}_t (r) = \int^\infty_0 \e^{-(\lambda^2 + \rho^2)t}
\phi_\lambda^{(\alpha, \beta)} (r) |c_{\alpha, \beta} (\lambda)|^{-2}
\d\lambda, \label{heat-kernel-jacobi-2}
\end{equation}
where
\begin{equation*}
c_{\alpha, \beta} (\lambda) = \f{2^{\rho - i \lambda} \Gamma (\alpha +
1) \Gamma (i \lambda)} {\Gamma (\f12 (i \lambda + \rho)) \Gamma (\f12 (i
\lambda + \alpha - \beta + 1))}.
\end{equation*}
We will see below that the type of $\alpha, \beta$ we are interested in
satisfy the above restriction of $\alpha, \beta$ and hence the inversion
formula is valid there.

We have the following sharp estimate of $h^{(\alpha, \beta)}_t$,
due to Anker {\it et~al} \cite{ADY}.

\setcounter{theo}{0}
\begin{theor}[\!] Let $\alpha, \beta\in \Z$ and $\alpha \ge \beta \ge -
\f12$. Then for $t > 0${\rm ,}
\begin{equation}
h^{(\alpha, \beta)}_t(r) \asymp t^{-3/2} \e^{-\rho^2 t} (1 + r) \left(1
+ \f{1 + r}{t} \right)^{\alpha - \f12} \e^{-\rho r} \e^{-r^2/4t}.
\label{heat-estimate}
\end{equation}
\label{ADY-estimate}
\end{theor}

$\left.\right.$\vspace{-1.5pc}

\noindent{\it Back to symmetric space.}\ \ Let us now come back to the
symmetric space $X$ and assume that $\alpha = \f{m_\gamma + m_{2
\gamma}-1}{2}$ and $\beta =\f{m_{2\gamma} - 1}{2}$. Then $\rho = \alpha
+ \beta + 1$ is the same as $\rho_0 = \f{m_\gamma + 2m_{2 \gamma}}{2}$.
In this case the operator ${\mathcal L}_{\alpha, \beta}$ coincides with
the radial part of the Laplace--Beltrami operator $\Delta$ of $X$. In the
previous section we have noticed that the usual heat kernel is
bi-invariant. Therefore, $\Delta h_t = {\mathcal L}_{\alpha, \beta} h_t$
and consequently, the heat kernel of ${\mathcal L}_{\alpha, \beta}$,
$h^{(\alpha, \beta)}_t(r)$ is the same as the heat kernel $h_t(a_r)$ of
$\Delta$. Note also that in this case, $\mu (\lambda) = |c
(\lambda)|^{-2} = |c_{\alpha, \beta} (\lambda) |^{-2} \asymp (1 +
|\lambda|)^{2 \alpha + 1}$ by (\ref{estimate-plancherel-measure}).

In all rank 1 symmetric spaces except the real hyperbolic space, $\alpha$
and $\beta$ satisfy the restriction in Theorem~5.1 and hence have the
estimate (\ref{heat-estimate}). (A full chart of $\alpha, \beta$ etc.
for the rank 1 symmetric spaces is available for instance in \cite{FK}.)
On the other hand for real hyperbolic space ${\mathcal H}^n = SO (n,
1)/SO(n)$, where $\alpha = \f{n-2}{2}$ and $\beta = - \f 12$, Davies and
Mandouvalos \cite{DM} already has the estimate of the type
(\ref{heat-estimate}) for the heat kernel $h_t$.

For $\delta \in \what{K}_0, h^\delta_t = h^{(\alpha + p_\delta, \beta +
q_\delta)}_t$ denotes the heat kernel corresponding to the operator
${\mathcal L}_{\alpha + p_\delta, \beta + q_\delta}$. By definition of
the heat kernel $J_\delta (h^\delta_t) = J_{\alpha + p_\delta, \beta +
q_\delta} (h^\delta_t) = \e^{-(\lambda^2 + \rho_\delta^2)t}$, where
$\rho_\delta = \alpha + p_\delta + \beta + q_\delta + 1$. Notice that
$\alpha + p_\delta, \beta + q_\delta$ are just two other $\alpha,
\beta$. When both $p_\delta, q_\delta$ are non-negative integers and
$\alpha, \beta$ satisfies the hypothesis of Theorem~5.1, $\alpha +
p_\delta, \beta + q_\delta$ will also satisfy the hypothesis and
consequently $h^\delta_t$ will again have the same estimate with
$\rho_\delta$ replacing $\rho$. There are however two exceptional cases.

\begin{enumerate}
\renewcommand\labelenumi{(\Alph{enumi})}
\leftskip .3pc
\item As noted earlier in the real hyperbolic spaces ${\mathcal H}^n$,
$\alpha, \beta$ do not satisfy the restriction in Theorem~5.1 and from
\cite{DM} we get the estimate (\ref{heat-estimate}) only for the
bi-invariant heat kernel $h_t$. To get a similar estimate for
$h^\delta_t$, where $\delta$ is not the trivial representation, we point
out that in $SO (n, 1), q_\delta = 0$ for all $\delta$ and that for
$(\alpha, \beta)$ of ${\mathcal H}^n$, $(\alpha + p_\delta, \beta)$ are
the corresponding parameters for the higher dimensional hyperbolic space
${\mathcal H}^{n + 2p_\delta}$ (see \cite{FK}). That is, $h_t^\delta$ of
${\mathcal H}^n$ is the same as $h_t$ of ${\mathcal H}^{n + 2
p_\delta}$. Thus we can again appeal to the result of \cite{DM} to show
that $h^\delta_t$ satisfies the estimate (\ref{heat-estimate})
substituting $\alpha$ by $\alpha + p_\delta$.

\item In the complex hyperbolic space $SU (n, 1)/S(U (n) \times U(1)),
p_\delta \ge 0, p_\delta > |q_\delta|$, but $q_\delta$ can be negative.
(In all other rank 1 symmetric spaces both $p_\delta$ and $q_\delta$
are non-negative.) To handle this exception, we shall denote by
$\tilde{\delta}$ the representation in $\what{K}_0$ which corresponds to
$(p_\delta, |q_\delta|)$. That is $|q_\delta| = q_{\tilde{\delta}}$. To
have a unform approach, in all symmetric spaces, instead of $(p_\delta,
q_\delta)$, we will consider $(p_\delta, |q_\delta|)$ and deal with
$h^{\tilde{\delta}}_t$, which will clearly satisfy
(\ref{heat-estimate}). Note that except for $SU (n, 1)/S(U(n) \times
U(1)), h^{\tilde{\delta}}_t$ is merely $h^{\delta}_t$ as $|q_\delta| =
q_\delta$.
\end{enumerate}

It is also clear that all $\alpha, \beta$ as well as $\alpha + p_\delta,
\beta + |q_\delta|$ we are concerned with satisfy the restriction for
the validity of the inversion formula (\ref{heat-kernel-jacobi-2}) (see
\cite{FK} p.~265 for $\alpha, \beta$'s and \cite{JW,Kos} for $p_\delta,
q_\delta$'s).

We proceed to find the heat kernel of the operator ${\mathcal
L}_{\alpha, \beta, p_\delta, q_\delta}$. The following expansion relates
the generalized spherical functions with the Jacobi functions:
\begin{align}
\Phi^1_{\lambda, \delta} (a_r) &= Q_\delta (\lambda) (\alpha +
1)_{p_\delta}^{-1} (\sinh r)^{p_\delta} (\cosh r)^{q_\delta}
\phi_\lambda^{\alpha + p_\delta, \beta + q_\delta} (r)\nonumber\\[.2pc]
&= Q_\delta (\lambda) (\alpha + 1)_{p_\delta}^{-1} \phi_{\lambda,
p_\delta, q_\delta}^{\alpha, \beta}(r),\label{Phi-delta}
\end{align}
where $x = ka_rK, Q_\delta (\lambda) = (\f12 (\alpha + \beta + 1 + i
\lambda))_{\f{p_\delta + q_\delta}{2}} (\f12 (\alpha - \beta + 1 + i
\lambda))_{\f{p_\delta - q_\delta}{2}})$ is the Kostant polynomial and
$(z)_m = \Gamma (z + m)/\Gamma(z)$. This relation is due to Helgason
(see \cite{He3}). We have used the parametrization given in \cite{B}. It
follows in particular from (\ref{Phi-delta}) that
$\phi_\lambda^{(\alpha, \beta)}(r)$ coincides with the elementary
spherical function $\phi_\lambda (a_r)$ and that $\Phi^1_{\lambda,
\delta}$ is an eigenfunction of ${\mathcal L}_{\alpha, \beta, p_\delta,
q_\delta}$ with the eigenvalue $-(\lambda^2 + \rho_0^2)$.

For a function $f$ on $X$ with a suitable decay, the $j$th
component of the {\em $\delta$-spherical transform} of $f$ is
\begin{equation*}
\what{f}_{\delta, j} (\lambda) = \int_{X} f(x) \Phi^j_{-\lambda, \delta}
(x) \d x
\end{equation*}
(see \cite{B,He2}). It is easy to see that $\what{f}_{\delta, j} (\bar
\lambda) = \what{f}_u (\lambda)$ defined in the previous section when $u
= Y_{\delta, j}$.

We define a function $H^\delta (x, t) = H^\delta_t(x)$ on $X \times
R^+$, through its $\delta$-spherical transform as
$\what{H^\delta}_{t\,\delta, j} (\lambda) = Q_\delta (\lambda)
\e^{-(\lambda^2 + \rho^2_0) t}$, for $1\le j \le d_\delta$ and
$\what{H^{\delta}}_{t\,\delta',j'}\equiv 0$ for $1\le j' \le
d_{\delta'}$, when $\delta'\neq \delta$. Clearly for each fixed $t > 0$,
$H^\delta_t$ is a function on $X$ of left $K$-type $\delta$. Since
$\Phi_{\lambda, \delta}^j$ is an eigenfunction of $\Delta$ with
eigenvalue $-(\lambda^2 + \rho^2_0)$, it is easy to see that
$H^\delta_t$ is a solution of type $\delta$ of the heat equation,
$\Delta f = \f{\partial}{\partial t}f$. Note that the Fourier transform
of $H^\delta_t$ has only the generic zeroes of the Fourier transforms of
functions of left-type $\delta$. In that way it is the basic solution
(of type $\delta$) of the heat equation and may be viewed as a
generalization of the bi-invariant heat kernel to arbitrary $K$-types.
We shall see that $H^\delta_t(x)$ also satisfies an estimate similar to
that in Theorem \ref{ADY-estimate}.

Let us recall that $Y_{\delta, j}$'s ($\delta \in \what{K}_0, 1 \le j\le
d_\delta$) are spherical harmonics. In analogy with (\ref{F_m-defn}),
for a nice measurable function $f$ on $X$ we can define, $F_{\delta, j}
(f) (\lambda) = Q_\delta (-\lambda)^{-1} \int_{K/M} \tilde{f} (\lambda,
kM) Y_{\delta, j} (kM) \d k$ and $\tilde{F}_{\delta, j} (f) (\lambda) =
Q_\delta (-\lambda) F_{\delta, j} (f) (\lambda)$.

\begin{rema}{\rm 
Looking back at the results of \S4 we see that explicit knowledge about
the matrix coefficients of the principal series representations may lead
to a refinement of the characterization given in Theorems~4.1 and 4.3.
For instance here when $G$ is of real rank 1, if $\deg Q_\delta \ge
P_{u,b}$ where $u$ transform according to $\delta \in \what{K}$ and
$Q_\delta$ is the Kostant polynomial defined above then $f_{u,0} = 0$.
Therefore if $\deg Q_\delta > \f{2}{p'} + \f{k}{p}$ or $\deg Q_\delta
\ge \f{l - 2 \alpha - 1}{q}$, then $f_{\delta} = 0$. Unfortunately no
exact description of the properties of the matrix coefficients of the
representations is available in general rank, because such a description
needs an exhaustive understanding of the subquotients of the principal
series representations.}
\end{rema}

With this preparation we are now in a position to state our last result:

\begin{theor}[\!]\label{rankone}
Let $f$ be a measurable function on $X$ and $p, q\in [1,\infty)$. Suppose
for some $k, l \in \R^+${\rm ,}
\begin{equation}
\int_{X} \f{|f(x) h_t(x)^{-1} \Xi (x)^{\f{2}{p}}|^p}{(1 + \sigma(x))^k}
\d x < \iy \label{onfunction-rankone}
\end{equation}
and for every $\delta\in \what{K}_0$ and $0\le j \le d_\delta${\rm ,}
\begin{equation}
\int_{\mathfrak a^*} \f{|F_{\delta,j} (f) (\lambda) \e^{t |\lambda|^2}|^q}{(1 + |\lambda|)^l}
|c(\lambda)|^{-2} \d\lambda < \iy.\label{onft-rankone}
\end{equation}
If also $l \le q + 2 \alpha + 2${\rm ,} then $f = \Sigma_{\delta \in
\what{K}_0 (\f{k-1}{p})} f_{\delta}$ is left $K$\!-finite where $f_\delta$
is the left $\delta$-isotypic component of $f$ and

\begin{enumerate}
\renewcommand\labelenumi{\rm (\alph{enumi})}
\leftskip .1pc
\item if for $\delta\in \what{K}_0 (\f{k-1}{p})$ and $1 \le j \le
d_\delta, F_{\delta,j} (f) (i \rho_0) = 1$ or $F_{\delta, j} (f) (-i
\rho_0) = 1$ according as $Q_\delta (i \rho_0) \neq 0$ or $Q_\delta (-i
\rho_0) \neq 0${\rm ,} then
\begin{equation*}
\hskip -.5cm f(a_r) = \Sigma_{\delta \in \what{K}_0 (\f{k-1}{p})}
H^\delta_t (a_r).
\end{equation*}

\item if $k \leq p + 1$ and $\int_{X} f(x)\d x = 1${\rm ,} then $f =
H^0_t = h_t$.
\end{enumerate}
\end{theor}

\begin{proof}
Let $f_{\delta, j} (a_r) = \int_K f(ka_r) Y_{\delta, j} (kM) \d k$. Then
$f_{\delta, j}$ satisfies (\ref{onfunction-rankone}).

From the definition of the Helgason Fourier transform
(\ref{Helgasontransform}) we have
\begin{equation*}
F_{\delta, j} (f) (\lambda) = Q_\delta (\lambda)^{-1} \int_{X}
\int_{K/M} f(x) \e^{(-i \lambda + \rho) A(x, b)} Y_{\delta, j} (b) \d
b\, \d x.
\end{equation*}
Let $x = ka_rK$. Using (\ref{helgason-to-jacobi}) and then changing over
to polar coordinates we get,
\begin{align*}
F_{\delta, j} (f) (\lambda) &= Q_\delta (\lambda)^{-1} \int_{X} f(ka_r)
Y_{\delta, j} (kM) \Phi^1_{-\lambda, \delta} (a_r)\d x\\[.2pc]
&= Q_\delta (\lambda)^{-1} \int_0^\infty f_{\delta, j} (a_r)
\Phi^1_{-\lambda, \delta} (a_r) \Delta_{\alpha, \beta} (r) \d r.
\end{align*}
Using (\ref{Phi-delta}) and (\ref{jacobi-delta}), we further have,
\begin{equation}
F_{\delta, j} (f) (\lambda) = \f{4^{p_\delta + q_\delta}}{(\alpha +
\beta)_{p_\delta}} \int_0^\infty f^\delta_{\delta, j} (a_r)
\phi_{-\lambda}^{(\alpha + p_\delta, \beta + q_\delta)}(r)
\Delta_{\alpha + p_\delta, \beta + q_\delta} (r) \d r, \label{F-delta-j}
\end{equation}
where
\begin{equation}
f^{{\delta'}}_{\delta,j} (a_r) = f_{\delta,j} (a_r) (\sinh
r)^{-p_{\delta '}} (\cosh r)^{-q_{\delta'}}. \label{f-tilde}
\end{equation}

Note that in the lone case (namely in $G = SU(n,1)$) where $q_\delta$
can be negative and hence $|q_\delta| \neq q_\delta, \beta$ is zero.
Therefore in view of the discussion (B) above we will rewrite
(\ref{F-delta-j}) using (\ref{jacobi-beta-to-minus-beta}) as
\begin{equation}
F_{\delta, j} (f) (\lambda) = C_\delta \int_0^\infty
f^{\tilde{\delta}}_{\delta, j} (a_r) \phi_{-\lambda}^{(\alpha +
p_\delta, \beta + |q_\delta|)} (r) \Delta_{\alpha + p_\delta, \beta +
|q_\delta|} (r) \d r.\label{F-delta-j2}
\end{equation}

We can now follow exactly the same steps as of Theorem~\ref{main-result}
and use (\ref{onfunction-rankone}), (\ref{estimate-jacobi-function}),
(\ref{jacobi-delta}) and finally appeal to Lemma~\ref{lemma-2} to show
that $F_{\delta, j} (f) (\lambda) = C^t_{\delta, j} \cdot
\e^{-t\lambda^2}$.

But by (\ref{F-delta-j2}), $F_{\delta, j} (f) (\lambda)$ is the Jacobi
transform of type $(\alpha + p_\delta, \beta + |q_\delta|)$ of
$f^{\tilde{\delta}}_{\delta, j}$. Therefore from
(\ref{heat-kernel-jacobi-2}) and (\ref{f-tilde}),
\begin{equation}
f_{\delta, j} (a_r) = C_{\delta, j}(\sinh r)^{p_\delta}(\cosh
r)^{|q_\delta|}h^{\tilde{\delta}}_t(r).
\end{equation}

On the other hand from (\ref{Helgasontransform}) and (\ref{Phi^j})
we have
\begin{align*}
\tilde{F}_{\delta, j} (f) (\lambda) &= \int_{K/M} \tilde{f} (\lambda,
kM) Y_{\delta, j} (kM) \d k\\[.2pc]
&= \int_{X} \int_{K/M} f(x) \e^{(-i \lambda + \rho) A (x, b)} Y_{\delta,
j}(b) \d b\, \d x\\[.2pc]
&= \int_{X} f(x) \Phi^j_{-\lambda, \delta} (x) \d x\\[.2pc]
&= \what{f}_{\delta, j}(\lambda).
\end{align*}
Therefore, $F_{\delta, j} (f) (\lambda) = C_{\delta, j} \cdot \e^{-t
\lambda^2}$ for $\lambda \in \R$ implies that $f_{\delta, j} (x) =
C^t_{\delta, j} H^\delta_{t, j} (x)$. That is,
\begin{equation}
f_{\delta, j} (a_r) = C^t_{\delta, j} (\sinh r)^{p_\delta} (\cosh
r)^{|q_\delta|} h^{\tilde{\delta}}_t (r) = C^t_{\delta, j} H^\delta_{t,
j} (a_r). \label{relate-jacobi-group-heatkernel}
\end{equation}
The apparent contradiction that arises on taking the Fourier transform
of the sides can be resolved by noting that $Q_\delta =
Q_{\tilde{\delta}}$ because when $|q_\delta| \neq q_\delta$ then $\beta
= 0$.

Now from Theorem~\ref{ADY-estimate} and the subsequent discussions
(A) and (B), it follows that
\begin{align}
H^\delta_t(a_r) &\asymp C(\delta, t) (1 + r) \left(1 + \f{1 + r}{t}
\right)^{\alpha + p_\delta - \f12} \e^{-(\rho_0 r +
\f{r^2}{4t})}\nonumber\\[.2pc]
&\asymp C(\delta, t) \left(1 + \f{1 + r}{t} \right)^{p_\delta} h_t(r).
\label{delta-estimate}
\end{align}
Notice that we have switched from $q_\delta$ to $|q_\delta|$ to fulfill
the requirement of Theorem~\ref{ADY-estimate}.

If $\delta \in \what{K}_0 \setminus \what{K}_0 (\f{k-1}{p})$, then
$p_\delta \ge \f{k-1}{p}$ and hence $f_{\delta, j} \equiv 0$ for $j = 1,
\dots, d_\delta$ by (\ref{onfunction-rankone}) and
(\ref{delta-estimate}).

Now (a) follows by taking Fourier transforms of the two sides of
$f_{\delta, j} (a_r) = C^t_{\delta, j} H^\delta_t (a_r)$ and putting the
initial condition.

If $k\leq p + 1$, then $p_\delta \ge \f{k-1}{p}$ whenever $\delta$ is
non-trivial as $p_\delta > |q_\delta|$. Therefore, $f = C_{t}h_t$. In
particular $f$ is a bi-invariant function. Note that $\int_{X} f(x) \d x
= \what{f} (i\rho_0)_{0,0}$, where $\what{f}(\cdot)_{0, 0}$ is the
spherical Fourier transform of $f$. Using the condition $\int_{X} f(x)
\d x = \what{f}(i\rho_0)_{0,0} = 1$, we get $f = h_t$. Thus (b) is
proved.\hfill $\Box$
\end{proof}

\section*{Acknowledgement}

We thank S~Helgason for sending us his survey paper \cite{HeNew} where
the preprint version of this paper is cited.

\end{document}